\documentclass [11pt]{article} 
  \usepackage{amsmath}

\newtheorem{theorem}{Theorem}[section]
\newtheorem{lemma}{Lemma}[section]
\newtheorem{corollary}{Corollary}[section]
\newtheorem{proposition}{Proposition}[section]
\newtheorem{remark}{Remark}[section]
\newcommand{\eqnsection}{
   \renewcommand{\theequation}{\thesection.\arabic{equation}}
   \makeatletter
   \csname @addtoreset\endcsname{equation}{section}
   \makeatother}
 
\def \ov{\overline}

\def \be{\begin{equation}}
\def \ee{\end{equation}}
\def \bt{\begin{theorem}} 
\def \et{\end{theorem}}
\def \bl{\begin{lemma}} 
\def \el{\end{lemma}}
\def \bea{\begin{eqnarray}} 
\def \eea{\end{eqnarray}}
\def \bas{\begin{eqnarray*}}
\def \eas{\end{eqnarray*}}

\def \al{\alpha}
\def \bb{\beta}
\def \ga{\gamma}
\def \Ga{\Gamma}
\def \de{\delta}

\def \ep{\epsilon}

\def \la{\lambda}

\def \om{\omega}
\def \Om{\Omega}

\def \si{\sigma}

\def \th{\theta}

\def \ff{\infty}
\def \wh{\widehat}
\def \wt{\widetilde}

\def\stl{\stackrel{\LL}{=}}
\def \cd{\,\cdot\,}

\def \cd{\,\cdot\,}

\def \AA{{\cal A}}
\def \BB{{\cal B}}

\def \FF{{\cal F}}
\def \GG{{\cal G}}

\def \II{{\cal I}}

\def \LL{{\cal L}}

\def \RR{{\cal R}}
\def \SS{{\cal S}}

\def \({\left(}
\def \){\right)}

\def \lc{\left\{}
\def \rc{\right\}}

\def \nn{\nonumber}
\def \Proof{\noindent{\bf Proof $\,$ }}

\def \bc{\begin{center} }
\def \ec{\end{center} }
\def \bs{\begin{slide} }
\def \es{\end{slide} }

\def\square{{\vcenter{\vbox{\hrule height.3pt
        \hbox{\vrule width.3pt height5pt \kern5pt
           \vrule width.3pt}
        \hrule height.3pt}}}}
\def\qed{{\hfill $\square$ \bigskip}}

\eqnsection
\begin{document}

\title{  Permanental processes  }

    \author{ Hana Kogan\,\,Michael B. Marcus\,\, Jay Rosen \thanks{Research of   the 
authors was supported by  grants from the National Science Foundation and of the second and third listed authors by grants from 
PSCCUNY.}}
\maketitle

\begin{abstract}    This is a survey of results about permanental processes, real valued positive processes which  are a generalization of squares of  Gaussian processes. In a certain sense the symmetric positive definite function that determines a Gaussian process is replaced by a function that is not   necessarily symmetric nor positive definite, but that nevertheless determines a   stochastic process.   This is a new avenue of research with very many   open problems. \end{abstract}

\bibliographystyle{amsplain}

   \section{Introduction}   \label{sec-1}
  
   Let $T$ be an index set and    $G=\{G(x);  x\in T\}$ be a mean zero Gaussian process with covariance  $U=\{U(x,y); x,y\in T\}$.  The covariance $U$ completely determines $G$. If instead of $G$ we decided to consider its square, i.e.,  $G^{2}:=\{G^{2}(x); x\in T\}$, we would not find this difficult, since we could simply obtain information about $G^{2}$ from our extensive understanding of the properties of $G$. However, let us suppose that we were completely unaware of the theory of Gaussian processes. We then might begin our study of $G^{2}$ by considering its finite joint distributions, as given by their Laplace transforms, i.e., 
    \begin{equation}
   E\(\exp\(-\frac{1}{2}\sum_{i=1}^{n}\al_{i}G^{2} (x_{i})\)\)=\frac{1}{|I+\al \Ga |^{1/2}}\label{1.8q} 
      \end{equation}
       for all $x_{1},\ldots,x_{n}$ in $T$,     where $I$ is the $n\times n$ identity matrix, $\al$ is the diagonal matrix with $\(\al_{i,i}=\al_{i}\)$,  $\al_{i}\in R_{+}$ and $\Ga=\{\Ga(x_{i}, x_{j})\}$ is an $n\times n$ matrix, that is symmetric and positive definite. (Of course, in truth, we do know the theory of Gaussian processes so we know that \{$\Ga(x_{i}, x_{j})\}$ is actually the covariance of the mean zero Gaussian vector, $(G(x_{1}), \ldots,G(x_{n}))$. That is, $\Ga=\{U(x_{i}x_{j})\}_{i,j=1}^{n}$.) Nevertheless let us continue to pretend that we do not know   the theory of Gaussian processes. A natural problem for us to consider is to  find necessary and sufficient conditions for $G^{2}$ to have continuous paths. Considering that the only unknown on the right-hand side of (\ref{1.8q}) is the matrix $\Ga$, the answer would have to be given in terms of $\Ga$.
       
       Lets now suppose that we have obtained necessary and sufficient conditions for $G^{2}$ to have continuous paths. We might then ask whether we could relax the conditions that $\Ga$ is symmetric and positive definite in (\ref{1.8q}) and that the power of  the determinant is 1/2.
       For example we might ask when we can define a real valued positive stochastic process  $\th =\{\th_{x}, x\in T\} $ with finite joint distributions that satisfy 
 \begin{equation}
   E\(\exp\(-\frac{1}{2}\sum_{i=1}^{n}\al_{i}\th_{x_{i}}\)\)=\frac{1}{|I+\al \Gamma |^{\bb}}\label{1.8w}, 
      \end{equation}
 where   $\Ga=\{\Ga(x_{i}, x_{j})\}_{i,j=1}^{n}$ is an $n\times n$ matrix and $\bb>0$.  Of course, aside from dimension 1 where $\th_{x_{1}}$ has a gamma distribution, for $\bb>0$, it is not at all clear whether (\ref{1.8w}) makes sense, i.e.,  the right hand side may not be a Laplace transform. 
 
   In 1997,      D. Vere-Jones,  \cite{VJ},  initiated the study of random vectors  with Laplace transform  given by (\ref{1.8w}).  He also briefly considered stochastic processes for which (\ref{1.8w}) gives finite joint distributions. We call these random vectors and stochastic  processes $\bb$-permanental processes.

In \cite{VJ} Vere-Jones gives necessary and sufficient conditions in terms of  $\Ga$ and $\al$ for (\ref{1.8w}) to exist, i.e., that  the right-hand side is a Laplace transform. Unfortunately, at least to us, it seems impossible to verify them, except in a very important case that was observed by N. Eisenbaum and H. Kaspi \cite{EK}. This is the case when $\Ga$ is the potential density of a  transient Markov process on $T$, or, equivalently, when all the finite dimensional matrices $\Ga$ are $M$-matrices, or, equivalently, when all finite dimensional vectors $(\th_{x_{1}},\ldots,\th_{x_{n}})$ are infinitely divisible. None of these three mathematical concepts is generally  considered by    experts in Gaussian processes. One purpose of this review is to convince the reader that they should   be.

 Several researchers have been looking for something like permanental processes for  a long time. The Dynkin Isomorphism Theorem, (DIT)  relates the squares of Gaussian processes to the local times of Markov processes, which necessarily have symmetric 0-potential densities. An outstanding problem was to find a version of Dynkin's theorem when the  0-potential densities of the Markov processes  are not symmetric. This is no doubt what really inspired Eisenbaum and Kaspi. They give a generalization of the DIT, or to be more precise, of a version of the DIT given in \cite{EKMRZ}, in \cite[Corollary 3.5]{EK}. Of course, the generalization of the process that is Gaussian squares, is the permanental process. The DIT enables one to use results about Gaussian processes to obtain results about the local times of associated Markov processes.  This is what \cite{book} is devoted to. In \cite[Theorem 1.3]{perm}   Theorem \ref{theo-1.1a} in this paper is used to obtain a sufficient condition for the continuity of local times of Markov processes
 without requiring that they have  symmetric 0-potential densities. This result with a slightly more restrictive hypothesis was obtained in \cite{EK1} by a different method. We do not discuss Dynkin type isomorphism theorems further in this paper; our focus is on permanental processes. However, we advise any researcher who decides to work on permanental processes to study them. 
 
 Another fascinating property of some permanental processes is that they can be represented as  the local times of loop soups associated with   Markov processes. This was shown for permanental processes associated with symmetric  Markov processes in \cite{LeJ}, (see also \cite{LeJ1}), and extended to permanental processes associated with more general   Markov processes in \cite{perm}.

\medskip	In Section \ref{sec-2} we give Verre-Jones's definition of $\bb$ permanents and his necessary and sufficient condition for the existence of $\bb$-permanental vectors. In Section \ref{sec-3}, under the assumption that (\ref{1.8w}) defines a permanental process, we give  sufficient conditions for them to be bounded or continuous. In  Section  \ref{sec-4} we summarize some of the results in \cite{EK} that show that $\bb$-permanental   processes exist if their kernels are the zero potentials of transient Markov processes. (In which case the permanental processes are infinitely divisible.) L\'evy processes killed the first time they hit zero, or at the end of and independent exponential time, have zero potentials that are kernels of permanental processes. In Section \ref{sec-pfbm} we obtain these potentials of L\'evy processes and thus get an idea of what the kernels of permanental processes, that are not the squares of Gaussian processes, look like. To be  specific, in Section \ref{subsec-5.1}, we give the kernels of permanental processes that are  generalizations of the squares of fractional Brownian motion. The study of permanental processes is just beginning. In Section \ref{sec-6} we list a few problems that interest us.

This paper is based on Eisenbaum and Kaspi's seminal paper \cite{EK} and on \cite{perm}. It is primarily a survey except that Remark \ref{rem-3.1}, Lemmas \ref{lem-4.2}--\ref{cor-1.1} and Theorem \ref{theo-4.2} are new results. Also the   process {$ \bf FBMQ^{\al ,\bb}$} appears for the first time in this paper.

\section{Existence of permanental processes} \label{sec-2}

The square of a Gaussian process is a permanental process so there are plenty of examples. But we are interested here primarily  in those that are not  squares of Gaussian processes.  To begin it is important to note that the right-hand side   of (\ref{1.8w}) is not unique with respect to $\Ga$. If $D$ is any diagonal matrix with non-zero entries we have
    \begin{equation}
   |I+\al \Ga |=|I+\al D^{-1}\Ga D |=|I+\al D^{-1}\Ga^{T} D |.
   \end{equation}
 For a very large class of   irreducible matrices $\Ga$, it is known that these are the only sources of non-uniqueness; see \cite{L}.
On the other hand, in certain extreme cases, for example, if $\Ga_{1}$ and $\Ga_{2}$ are $n\times n$ matrices with the same diagonal elements and all zeros below the diagonal, then $|I+\al \Gamma_{1} |=|I+\al \Gamma_{2} |$.  Therefore, a process that is the square of a Gaussian process can have a   non-symmetric kernel. Nevertheless, as we    show, there are many permanental processes that are not squares of Gaussian processes. 

The first step in the definition of permanental processes is the definition of the $\bb$-permanent of a $n\times n$ matrix $B=\{b_{i,j}\}_{i,j=1}^{n}$, which we denote by $|B|_{\bb}$. It is
\begin{equation}
   |B|_{\bb}=\sum_{\si}\bb^{m(\si)}b_{1,i_{1}}b_{2,i_{2}}\ldots b_{n,i_{n}},
   \end{equation}
where the summation is taken over all distinct permutations
\begin{equation}
   \si= {1,\,2,\ldots,n\choose i_{1}\,,\,i_{2},\ldots ,i_{n} }\label{2.2}
   \end{equation}
of the indices $1,2,\ldots,n$  and $m(\si)$ denotes the number of cycles into which $\si$ can be decomposed. 

The name $\bb$-permanent is taken from the fact that $|B|_{1}$ is the permanent of $B$.

\medskip	For $\bb>0$  we define what it means for  an  $n\times n$ matrix $B$ to   be $\bb$ positive definite. This is a rather complicated definition. To begin consider   all multi-indices ${\bf k}=(k_{1},\ldots,k_{n}) $, (where the $k_{i}$ are integers), and set $|{\bf k}|=k_{1}+\cdots+k_{n} $. For each ${\bf k}$, let  $B({\bf k})$ denote the ${|\bf k|}\times |{\bf k}|$ matrix   with entries
\begin{equation}
   \{B({\bf k})\}_{i,j}=B_{p_{i },p_{j }}\label{2.5}
   \end{equation}
where $p_{i}=1$ if $1\le i\le k_{1}$, and $p_{i}=l$ if $k_{l-1}\le l\le k_{l}$, if $l=2,\dots,n$.
The matrix $B$ to said to   be $\bb$ positive definite if for all multi-indices $|{\bf k}|$
\begin{equation}
   |B({\bf k})|_{\bb}\ge0.\label{2.5s}
   \end{equation}
 Clearly, to verify whether $B$ is $\bb$ positive definite one must check (\ref{2.5s}) for an infinite number of matrices.
 
\medskip	In \cite[Proposition 4.5]{VJ} Verre-Jones gives necessary and sufficient conditions for (\ref{1.8w}) to be the Laplace transform of the vector $(\th_{x_{1}},\ldots,\th_{n})$ in terms of the  modified resolvent matrix 
\begin{equation}
   \Ga_{r}:= \Ga(I+r\Ga)^{-1}\label{2.6}
   \end{equation}
where  $r\ge 0$, and $\Ga$ is the matrix in (\ref{1.8w}).   Note that when $\Ga^{-1}$ exists, we have   $   \Ga_{r}=  (\Ga^{-1}+r I)^{-1}$.

 The next proposition is Proposition 4.5 in \cite{VJ}.

\begin{proposition} \label{prop-2.1}For (\ref{1.8w}) to represent the Laplace transform of a non-negative random vector it is necessary and sufficient that for all  $r\ge 0$
\begin{itemize}
\item[(i)] $\Ga_{r}$ exists and is $\bb$-positive definite.
\item[(ii)] $ \det(I+r\Ga)>0$.
\end{itemize}
\end{proposition}

(Item {\em (ii)} is equivalent to: All the real, non-zero, eigenvalues of $\Ga$ are positive.)

\medskip	Verifying Proposition \ref{prop-2.1} {\em (i)}, requires verifying (\ref{2.5}) for an   infinite number of matrices,  which seems impossible unless   all the entries of the matrix $\Ga_{r}$ are greater than or equal to zero.  In \cite[Theorem 3.1]{EK} Eisenbaum and Kaspi point out that this is is the case for an important class of kernels $\Ga$ and that Proposition \ref{prop-2.1}  {\em (ii)}, also holds for these kernels. We take this up in Section  \ref{sec-4}. In the next section we assume that a permanental process exists and give sufficient conditions for it to be continuous.
 
\section{Continuity of permanental processes}\label{sec-3}

   A key observation that allows us to obtain a `best possible' sufficient condition for the continuity of permanental processes is that the bivariate marginals of a $1/2$-permanental process are squares of  bivariate normal random variables.  
We proceed to explain this.

For $n=2$, (\ref{1.8q}) takes the form  
  \bea
  && E\(\exp\(-\frac{1}{2}\(\al_{1} \th_{x }+\al_{2} \th_{y}\)\)\)\label{1.8a}\\
  &&\qquad=\frac{1}{|I+\al \Gamma |^{1/2}}\nn= \(1+\al_{1}\Ga(x , x )+\al_{2}\Ga(y, y)\right.\\
& &   \nn \qquad\qquad\left.+\al_{1}\al_{2}\(\Ga(x , x )\Ga(y, y)-\Ga(x , y)\Ga(y, x )\)\)^{-1/2}.\eea
An obvious necessary condition for this to exist is that $|I+\al \Gamma |>0$. Therefore, if we set  $\al_{2}=0$ in (\ref{1.8a})   and take $\al_{1}$  sufficiently large we see that for any   $x\in T$
 \begin{equation}
 \Ga(x , x )\geq 0\label{vj.1},
 \end{equation}
and taking $\al_{1}=\al_{2}$  sufficiently large,   that
    \be 
\Ga(x, x)\Ga(y, y)-\Ga(x, y)\Ga(y, x)\geq 0. \label{3.2}
 \ee
 In addition, by \cite[p. 135, last line]{VJ},  a necessary condition for the existence of any $\bb$-permanental process is that for any pair $x ,y\in T$  
  \begin{equation}
 \Ga(x , y ) \Ga(y , x )\geq 0.\label{vj.2}
 \end{equation}
  This allows us to define
     \begin{equation}
   d(x,y) =4\sqrt{ 2/3}\(\Ga(x,x)+\Ga(y,y)-2\(\Ga(x , y)\Ga(y, x )\)^{1/2}\)^{1/2}.\label{1.2w}
   \end{equation}
 It follows from (\ref{vj.1})--(\ref{vj.2}) that for any pair $x , y\in T,$   the matrix 
 \be\begin{bmatrix}
\Ga(x, x)& \(\Ga(x, y)\Ga(y, x)\)^{1/2}\\ 
\(\Ga(x, y)\Ga(y, x)\)^{1/2}&\Ga(y, y) \label{matrix}
\end{bmatrix}\ee
  is positive definite. Therefore we can construct a mean zero Gaussian vector $\{G (x),G (y)\}$ with covariance matrix 
 \begin{equation}
    E\(G(x) G(y)\) = \(\Ga(x, y)\Ga(y, x)\)^{1/2}.\label{1.11a}
   \end{equation}
  Note that  
     \begin{equation}
   \(E(G(x)-G(y))^{2}\)^{1/2}=\frac{\sqrt{3/2}}{4}d(x,y),.\label{dd}   \end{equation}
 
By (\ref{1.8a}) the Laplace transform of   $\{\th_{x },\th_{y}\}$  is   the same as 
  the Laplace transform of   $\{G^{2}(x ),G^{2}(y)\}$. This gives the following critical lemma:

 \begin{lemma}\label{lem-1}    Suppose  that  $\th:=\{\th_{x}, x\in T\} $ is a  1/2-permanental process
with kernel $\Ga$ as given in (\ref{1.8q}). Then for any pair $x , y$, 
\begin{equation}
   \{\th_{x},\th_{y}\}\stl  \{G^{2}(x),G^{2}(y)\}\label{fund}
   \end{equation}
where $\{G (x),G (y)\}$ is a mean zero Gaussian random variable with covariance matrix given by   (\ref{1.11a}).
 \end{lemma}
 
 For  $p\geq 1$, let $ 
\psi_p(x)= 
\exp(x^{p}) -1
$   and   $L^{\psi_p}(\Om,\FF,P)$ denote the set ofÊ random 
variables 
$\xi:\Om\to R^{1}$ such that $E\psi_p\( |\xi|/c \)<\ff$ for some $c>0$.  
$L^{\psi_p}(\Om,\FF,P)$ is a Banach space with norm given by
\be
\|\xi\|_{\psi_p}=\inf\left\{c>0:E\psi_p\( |\xi|/c\)\le1\right\}.\label{1.4}
\ee
We shall only be concerned with the cases $p=1$ and 2. 

It follows from Lemma \ref{lem-1} that 
 \begin{equation}
   \|\th_{x}-\th_{y}\|_{\psi_{1}}=\|G^{2}(x)-G^{2}(y)\|_{\psi_{1}}\le C\sup_{x\in T}\Ga^{1/2}(x,x)d(x,y)
   \end{equation}
  This is a good start but not the result we want   since we   know that the continuity condition when the permanental process is the square of a Gaussian process is the sufficient condition for continuity of the Gaussian process itself. This necessitates working with the Orlicz space $ L^{\psi_2} $. 
  
  Unfortunately,  $\|G^{2}(x)-G^{2}(y)\|_{\psi_{2}}=\ff$. We get around this by using a truncation argument introduced by Martin Barlow to obtain a sufficient condition for the continuity of local times of L\'evy processes \cite{Bar,Bar84}. The following lemma is \cite[Lemma 3.3]{perm}.

  \begin{lemma} \label{lem-3.3q}
Let $\th:=\{\th_{x}, x\in T\} $ be a  $1/2$-permanental process with kernel $\Ga$.  Then for all $x,y\in T$ and $0<  \la<\ff$   
 
 \be 
  \Big \| {\th_{x}\wedge\la \over \la^{1/2}}-{\th_{y}\wedge\la \over \la^{1/2}} \Big\|_{\psi_2}\le 
  d(x,y). \label{3.9a}
   \ee
 \end{lemma}

 The proof of Theorem \ref{theo-1.1a} below uses (\ref{3.9a}) and standard ideas used to prove the continuity of Gaussian processes along with some of the ideas in \cite{Bar,Bar84}.

 \medskip	Let $(T,\rho)$   be aÊ separable metric or pseudometric space.  Let $B_{\rho}(t,u)$  denote the closed ball in $(T,\rho)$  with radius $u$ and center $t$. For any probability measure  $\mu$  on $(T,\rho)$ we define
\begin{equation}
  J_{T, \rho,\mu}( a) =\sup_{t\in T}\int_0^a\(\log\frac1{\mu(B_{\rho}(t,u))}\)^{1/2} \,du.\label{tau}
   \end{equation}
   We occasionally omit    some of the subscripts $T, \rho$ or $\mu$, if they are  clear from the context.

  In general, $  d (x,y)$ is not a metric or   pseudometric on $T$. Nevertheless, we can still define the sets $B_{  d}(s,u)=\{t\in T\,|\,   d (s,t)\leq u\}$.  
We can then define $  J_{T,  d,\mu}( a)$ as in (\ref{tau}), for any probability measure $\mu$Ê Ê on $\mathcal{B}(T, d)$, the $\si$-algebra generated by the sets $B_{  d}(s,u)$.

 \bt\label{theo-1.1a}Ê Ê  Let $\th= \{\th_x: x\in T\}$Ê be a $1/2$-permanental process with kernel $\Ga$  satisfying $\sup_{x\in T} \Ga(x,x)<\ff $.  Let $D$ denote the  d-diameter of $T$ and assume that $T$ is separable for $  d  $, and that there exists
a probability measure $\mu$ÊÊ Ê on $\mathcal{B}(T, d)$ such that  
\be
J_{  d }(D)<\ff\label{1.1va}.
\ee 
Then there exists a version $\th'=\{\th'_x,x\in T\}$ of $\th$ such that   for any $x_{0}\in T$  
\be
 \|\sup_{x\in T}\th' _{x}  \|_{\psi_{1}}\le   4\| \th' _{x_{0}}  \|_{\psi_{1}}\label{1.5m} +C \(\sup_{x\in T}\Ga(x,x)\) 
J^{2}_{  d}(D),   
 \ee 
 where $C$ is a     constant. 
%Then there exists a version $\th'=\{\th'_x,x\in T\}$ of $\th$ which is bounded almost surely.
 
If
\begin{equation}
   \lim_{\de\to 0}J_{  d}(\de)=0,\label{1.8jv}
   \end{equation}
there exists a version $\th'=\{\th'_x,x\in T\}$ of $\th$ such that
 \be     \lim_{\de\to 0}\sup_{\stackrel{s,t\in T}{  d(s,t)\le \de}} |\th' _{s}( \om)-\th' _{t}( \om)|\label{3.6.90hhv} =0,\qquad a.s.
\ee 

If (\ref{1.8jv}) holds and 
\begin{equation} 
 \lim_{\de\to 0}{J_{  d}(\de)\over \de}=\ff
  \label{2.18va},
   \end{equation}
   then     
   \begin{equation}
   \lim_{\de\to0}   \sup_{\stackrel{s,t\in T}{  d (s,t)\le \de}}\frac{|\th' _{  s} -\th' _{t}|}{J_{  d} (  d(s,t) /2) }\le 30\(\sup_{x\in T}\th'_x\)^{1/2}\quad a.s.\label{2.1sv}
   \end{equation}
   \et

 Continuity is not mentioned in Theorem \ref{theo-1.1a} because we do not know whether $d$ is a metric on $T$. If $  d (x,y)$ is   continuous on $T\times T$ and   there exists a probability measure  $\mu $ on $T$ such that 
(\ref{1.8jv}) holds then  there exists a version $\th'= \{\th'_x: x\in T\}$ of $\th$ that is continuous almost surely. In this case (\ref{2.1sv}) gives a bound for the uniform modulus of continuity. It s not unreasonable to assume that $  d (x,y)$ is   continuous on $T\times T$. We show in \cite[Lemma 3.2]{perm} that when $\th$ is continuous on $T$ almost surely,   $  d (x,y)$ is continuous on $T\times T$.

It is easy to see why the term $\sup_{x\in T}\(\th'_x\)^{1/2}$ appears on the left-hand side of (\ref{2.1sv}). Consider a Gaussian process $  \{G_x: x\in T\}$ with modulus of continuity $\om(d(s,t))$. (Suppose that $(T,d)$ is a metric space.) Then
  \begin{equation}
   \lim_{\de\to0}   \sup_{\stackrel{s,t\in T}{  d (s,t)\le \de}}\frac{|G^{2}_{  s} -G^{2} _{t}|}{\om (  d(s,t)  ) }\le   \lim_{\de\to0}   \sup_{\stackrel{s,t\in T}{  d (s,t)\le \de}}\frac{|G _{  s} -G _{t}|}{\om (  d(s,t)  ) } \(2\sup_{x\in T}G^{2}_x\)^{1/2}\quad a.s. 
   \end{equation}
(With mild regularity conditions the inequality can be replaced by equality.)   

Under the hypotheses of Theorem \ref{theo-1.1a} an essentially optimal  local modulus is also given for $\th= \{\th_x: x\in T\}$ in \cite[Theorem 4.2]{perm}.
   
 \medskip	We  say that a metric or pseudometric $d_{1}$ dominates $  d$ on $T$
  if 
    \begin{equation}
     d(x,y)\leq d_{1}(x,y), \hspace{.3 in}\forall x,y\in T.\label{dom}
   \end{equation}
   In the \cite[Section 5]{perm}  we give several  natural  metrics that dominate $  d$.    For example
   \begin{equation}
\frac{1}{4}\inf_{x\in T}\Ga^{1/2} (x,x) d(x,y)\le    \( E\(\th _{x}-\th _{y} \)^{2}\)^{1/2}\le \frac{3 }{4}\sup_{x\in T}\Ga^{1/2}(x,x) d(x,y).
   \end{equation}
 
  If a metric $d_{1}$ dominates $d$ and is such that  $(T, d_{1})$ is separable and has finite diameter $D$ then the results of Theorem \ref{theo-1.1a} hold with $d$ replaced by $d_{1}$ and the version $\th'$ is uniformly continuous on $(T,d_{1})$.  (See \cite[Corollary 1.2]{perm}.)

\begin{remark} \label{rem-3.1}{\rm 

	We expect that Theorem \ref{theo-1.1a} also holds for all $\bb$-permanental processes, $\bb>0$, but we can't show this. What we can show is that 
all the statements in  Theorem \ref{theo-1.1a} and in the discussion following it hold for all $\bb$-permanental processes, $\bb>0$, if $  J_{T, \rho,\mu}( a)$ in  (\ref{tau}) is replaced by
 \begin{equation}
  \II_{T, \rho,\mu}( a) =\sup_{t\in T}\int_0^a\(\log\frac1{\mu(B_{\rho}(t,u))}\)  \,du.\label{tau1}
   \end{equation}
 This follows from the next lemma and corollary used in conjunction with the proof of  \cite[Theorem 1.1]{perm}. }\end{remark}
 
 \begin{lemma} \label{lem-3.3}
Let $\th:=\{\th_{x}, x\in T\} $ be a  $\bb$-permanental process, $\bb\le 1/2$  with kernel $\Ga$.  Then for all $x,y\in T$ and $0<  \la<\ff$  
 \be 
  \Big \| {\th_{x} }-{\th_{y} } \Big\|_{\psi_1}\le 
  \frac{16(1-\bb)}{3}\(\sup_{x\in T} \Ga (x,x) \)d(x,y). \label{3.9ab}
   \ee
 \end{lemma}

 \Proof     Let  $\Ga(x,y)$ be the kernel of the 2-dimensional vector $\th=(\th_{x},\th_{y})$. Let $G=(G_{x},G_{y})$ be the 2-dimensional Gaussian vector with covariance given by (\ref{1.11a}). By (\ref{1.8a}) we can assume that the kernel of $\th$ is given by the matrix in (\ref{matrix}) which we denote by $\Ga$. It is well known that  
  \begin{equation}
   E\(\exp\(-\frac{1}{2}(\al_{x}\th_{x }+\al_{y}\th_{y} )\)\)=\frac{1}{|I+\al \Gamma |^{\bb}}\label{1.8ww}, 
      \end{equation}
is the Laplace transform of vector in $R^{2}_{+}$ for all $\bb>0$. (See Lemma \ref{lem-4.1} and Remark \ref{rem-4.1}.)  Denote this vector by $ (\th_{x,\bb},\th_{y,\bb})$. For $\bb<1/2$ we have
\begin{equation}
   (\th_{x,\bb}-\th_{y,\bb})+ (\th'_{x,1/2-\bb}-\th'_{y,1/2-\bb})\stl G^{2}_{x}-G^{2}_{y}
   \end{equation}
where $ (\th'_{x,1/2-\bb},\th'_{y,1/2-\bb})$ is independent of  $(\th_{x,\bb},\th_{y,\bb})$, since $ (G^{2}_{x},G^{2}_{y})$ is the vector obtained from (\ref{1.8ww}) when $\bb=1/2$. Let $E_{\th}E_{\th'}$ denote expectation on the product space of $ (\th_{x,\bb},\th_{y,\bb})\times  (\th'_{x,1/2-\bb},\th'_{y,1/2-\bb})$ and let    $\|\cd\|_{p}= \(E_{\th}E_{\th'}|\cd|^{p}\)^{1/p}$, $p\ge 1$. We have
\bea
 && \|(\th_{x,\bb}-\th_{y,\bb})+ (\th'_{x,1/2-\bb}-\th'_{y,1/2-\bb}) \|_{p}\label{3.26}\\
 &&\qquad\nn =  \( E_{\th}E_{\th'}| (\th_{x,\bb}-\th_{y,\bb})+ (\th'_{x,1/2-\bb}-\th'_{y,1/2-\bb})|^{p}\)^{1/p}\\
  &&\qquad\nn \ge \( E_{\th}| (\th_{x,\bb}-\th_{y,\bb})+E_{\th'} (\th'_{x,1/2-\bb}-\th'_{y,1/2-\bb})|^{p}\)^{1/p}\\
 &&\qquad\nn \ge   \| \th_{x,\bb}-\th_{y,\bb}+E  (\th'_{x,1/2-\bb}-\th'_{y,1/2-\bb})\|_{p} \\
  &&\qquad\nn \ge   \| \th_{x,\bb}-\th_{y,\bb}\|_{p}-|E_{\th'} (\th'_{x,1/2-\bb}-\th'_{y,1/2-\bb})|.
   \eea
Consequently
\begin{equation}
   \|  \th_{x,\bb}-\th_{y,\bb} \|_{p}\le   \| G^{2}_{x}-G^{2}_{y} \|_{p} +|E_{\th'} (\th'_{x,1/2-\bb}-\th'_{y,1/2-\bb})|.\label{3.27}
   \end{equation}
 Using (\ref{1.8ww}) we see that for any $\bb<1/2$
 \bea
  |E_{\th'} (\th'_{x,1/2-\bb}-\th'_{y,1/2-\bb})|&=&(1-2\bb) |\Ga(x,x)-\Ga(y,y)|\\
   & \le&\nn  (1-2\bb) \| G^{2}_{x}-G^{2}_{y} \|_{p},
   \eea
  since for integers $p\ge 1$,
\begin{equation}
   \(E|G^{2}_{x}-G^{2}_{y}|^{p}\)^{1/p}\ge   \( |E\(G^{2}_{x}-G^{2}_{y}\)|^{p}\)^{1/p} = |\Ga(x,x)-\Ga(y,y)|.  \end{equation}
   Therefore
   \begin{equation}
   E(\th_{x,\bb}-\th_{y,\bb})^{p}\le (2(1- \bb))^{p}E(G^{2}_{x}-G^{2}_{y})^{p},\label{3.29}
   \end{equation}
which implies that 
\begin{equation}
    \| {\th_{x} }-{\th_{y} }  \|_{\psi_1}\le 2(1- \bb)  \| G^{2}_{x}-G^{2}_{y} \|_{\psi_1}
   \end{equation}
which gives (\ref{3.9ab}).\qed

\begin{corollary} Let $\th:=\{\th_{x}, x\in T\} $ be a  $\bb$-permanental process, $\bb>1/2$,  with kernel $\Ga$. Then for some integer  $k\ge 1$  
 \be 
  \Big \| {\th_{x} }-{\th_{y} } \Big\|_{\psi_1}\le 
  \frac{16(k-\bb)}{3}\(\sup_{x\in T} \Ga (x,x) \)d(x,y). \label{3.31}
   \ee
 \end{corollary} 
   
\Proof   Set $\bb=k\bb'$ for some $\bb'\le 1/2$. Then since
\begin{equation}
     \| {\th_{x,\bb} }-{\th_{y,\bb} }  \|_{\psi_1}\le k     \| {\th_{x,\bb'} }-{\th_{y,\bb'} }  \|_{\psi_1},
   \end{equation} 
   (\ref{3.31}) follows from (\ref{3.9ab}).  \qed

  In general, given the existence of a $1/2$-permanental process with kernel $\Ga$, we do not know whether there exists any other $\bb$-permanental process with kernel $\Ga$. If one does exist for rational $\bb$,  the next theorem shows that it inherits  continuity  properties of the $1/2$-permanental process. 
  
\medskip	We use the following inequality in the proof:

 \begin{lemma}\cite[Theorem 1.1.5]{DG} \label{lem-3.4}Let $X_{1},\ldots,X_{n}$ be independent identically distributed random variables in a normed linear space. Then for all $1\le k\le n$ and all $t>0$
\begin{equation}
   P\(\max_{1\le k\le n} \|\sum_{j=1}^{k}X_{j}\| >t\)\le 
9 P\( \|\sum_{j=1}^{n}X_{j} \| >{t\over 30}\).
   \end{equation}
   In particular
\begin{equation}
   P\(  \| X_{1}\| >t\)\le 
9 P\( \|\sum_{j=1}^{n}X_{j} \| >{t\over 30}\).\label{3.34}
   \end{equation}

  \end{lemma}

 \begin{theorem} \label{theo-3.2} Let $\th= \{\th_x: x\in T\}$Ê be a $1/2$-permanental process with kernel $\Ga$ that has  continuous paths on $(T,d)$ almost surely.   Let $\th_{\bb}= \{\th_{x,\bb}: x\in T\}$ be a $\bb$-permanental process with kernel $\Ga$ for any rational $\bb>0$. Then $\th_{\bb}$
has  continuous paths on $(T,d)$ almost surely.   
 \end{theorem}

   \Proof      Let $T'$ be a finite subset of  $T$ and $\|\cd\|:=\sup_{\stackrel{d(x,y)\le \de}{x,y \in T'}}|\cd|$.  Suppose that  $\bb=p/q$, for integers $p$ and $q$. Then  
 \bea
   \|\sum_{j=1}^{q}\({\th_{x,\bb} }^{(j)}-{\th_{y,\bb} }^{(j)}\)\|&\stackrel{\mathcal{L}}{=}&   \|\sum_{j=1}^{2p}{\th_{x,1/2} }^{(j)}-{\th_{y,1/2} }^{(j)}\|.\label{3.36}
   \eea 
   where the notation $(j)$ denotes independent copies. By (\ref{3.34}) and (\ref{3.36}) 
   \bea
&&   P\(  \| {\th_{x,\bb} } -{\th_{y,\bb} }\| >t\) \label{3.37} \\
   & & \qquad\le\nn
9 P\(   \|\sum_{j=1}^{2p}{\th_{x,1/2} }^{(j)}-{\th_{y,1/2} }^{(j)}\| >{t \over 30}\) \\ 
 & &\nn\qquad\le
 18p P\( \|{\th_{x,1/2} }  -{\th_{y,1/2} } \| >{t\over 60p}\)\\
 & &\nn\qquad=
 18p P\(\sup_{\stackrel{d(x,y)\le \de}{x,y \in T }}|{\th_{x,1/2} }  -{\th_{y,1/2} } | >{t\over 60p}\).
   \eea
   Since  $\th_{x,1/2}$ has continuous paths it follows from (\ref{3.37}) that for any finite subset  $T'$ of $T$ and $\ep>0$ we can take $\de$ sufficiently small so that 
   \begin{equation}
    P\(  \| {\th_{x,\bb} } -{\th_{y,\bb} }\| >\ep\)\le \ep.
   \end{equation}
  This implies that $\th_{\bb}$
has  continuous paths on $(T,d)$ almost surely. \qed

\section{Examples of permanental vectors and processes} \label{sec-4}

 The matrix
$A$  is said to be an  
$M$ matrix   if
\begin{enumerate}
\item[(1)] $a_{ i,j}\leq 0$ for all $i\neq j$.
\item[(2)] $A$ is nonsingular and $A^{ -1}\geq 0$.
\end{enumerate}

\begin{lemma} \label{lem-4.1}\cite[Lemma 4.2]{EK} Let $\{\Ga_{ i,j}\}$  be a real positive non-singular $n\times n$ matrix. There exists a positive infinitely divisible random vector $(\th_{1},\ldots,\th_{n})$ such that  for every $(\al_{1},\ldots,\al_{n})$ in $R_{+}^{n}$
 \begin{equation}
   E\(\exp\(-\frac{1}{2}\sum_{i=1}^{n}\al_{i}\th_{ {i}}\)\)=\frac{1}{|I+\al \Gamma |^{1/2}}\label{4.1qq}
      \end{equation}
if and only if $\Ga^{-1}$ is an $M$-matrix.
 \end{lemma}

\begin{remark}\label{rem-4.1}{\rm One can replace the $1/2$ in (\ref{4.1qq}) by any $\bb>0$. To see this denote the   right hand side of (\ref{4.1qq}) by $F(\al )$. A necessary and sufficient condition for a function $F(\al )$  to be the Laplace transform of an infinitely divisible random variable is that 
$   \log F(\SS) $, where $\SS$ is a diagonal matrix with entries  $(1-s_{i})t$,  $t>0$ and $0<s_{i}\le 1$, $i=1,\ldots,n$, has a power series expansion in $(s_{1},\ldots,s_{n})$  with positive coefficients for all $t$ sufficiently large. (See, e.g., \cite[Lemma 13.2.2]{book}.) In (\ref{4.1qq})
 \begin{equation}
   \log F(\SS) = -(1/2)\log|I+\SS \Gamma |.
   \end{equation}
 Clearly, if this has a power series expansion in $(s_{1},\ldots,s_{n})$  with positive coefficients for all $t$ sufficiently large, it does if 1/2 is replaced by any $\bb>0$.
 }\end{remark}

\medskip	Let $X=\{X_{t}, t\in R_{+}\}$ be a transient Markov process with state space $E$ and potential density $u(x,y)$.

\bt \label{theo-4.1}\cite[Theorem 3.1]{EK} Let $u(x,y)$ be the potential density of a transient Markov process with state space $E$. Then for every $\bb>0$ there exists  a positive process $\{\th_{x},x\in E\}$, such that  for every $(\al_{1},\ldots,\al_{n})$ in $R_{+}^{n}$ and $(x_{1},\ldots,x_{n})\in E$,
 \begin{equation}
   E\(\exp\(-\frac{1}{2}\sum_{i=1}^{n}\al_{i}\th_{x_ {i}}\)\)=\frac{1}{|I+\al U |^{\bb}}\label{4.1}
      \end{equation}
where $U=\{u(x_{i},x_{j})\}_{i,j=1}^{n}$.
 \et

The proof of this theorem makes use of the fact that because $u(x,y)$ is the potential density of a transient Markov process, $U^{-1}$ is an $M$-matrix. It is a property of $M$-matrices that all its real eigenvalues are positive. Consequently all the real eigenvalues of $U$ are positive. In addition the matrix $rU_{r}$, see  (\ref{2.6}), is a resolvent matrix and hence it has only non-negative entries. In this case it is trivial to see that  $U_{r}$ is $\bb$-positive definite for all $\bb>0$. Therefore, Theorem \ref{theo-4.1} follows from Proposition \ref{prop-2.1}.

\medskip	When a permanental process is determined by the potential of a transient Markov process we refer it it as an associated process. (This is ambiguous  because the kernel of a permanental process is not unique.)   
So far the only permanental vectors and processes that we can show exist are associated processes and they have the additional property that they are infinitely divisible. 

There is another property of associated processes that enables us to simplify the statement of Theorem \ref{theo-1.1a}. In general we don't know whether the function $d$ in this theorem is a metric on $T$. We show in the next lemma that when the permanental process is an associated process $d$ is a metric on $T$. In this case we can add to the statement of Theorem \ref{theo-1.1a} that the version $\th'$ is continuous on $(T,d)$ and that (\ref{2.1sv}) gives a uniform modulus of continuity on $(T,d)$.

\begin{lemma}\label{lem-4.2} Let $\{\th_{x}, x\in E\}$ be a $\bb$-permanental process associated with a transient Markov process with potential $u(x,y)$. Then
     \begin{equation}
  \wt  d(x,y) = 4\sqrt{2/3}\(u(x,x)+u(y,y)-2\(u(x , y)u(y, x )\)^{1/2}\)^{1/2}.\label{4.2w}
  \ee
is a metric on $E$. \end{lemma}

The proof follows from the next two lemmas:

  \bl \label{lem-1.1}Let 
\be  \AA=\left (   \begin{array}{cccc}
    u& a& b\\ 
   a& v& c  \\ 
 b&c & w \\ 
  \end{array} \right ),\qquad \BB=\left (   \begin{array}{cccc}
    u& \al_{1}& \bb_{2}\\ 
  \al_{2}& v& \ga_{1} \\ 
 \bb_{1}&\ga_{2} & w \\ 
  \end{array} \right ).
\ee
where $\al_{1}\al_{2}=a^{2}$, $\bb_{1}\bb_{2}=b^{2}$ and $\ga_{1}\ga_{2}=c^{2}$. Suppose that $\AA\ge 0$. If $\BB^{-1}$ is an $M$ matrix then $\AA^{-1}$ is an $M$ matrix.
\el

\Proof   We first consider the case when $abc>0$. We have
\begin{equation}
   \BB^{-1}=\frac{1}{\det \BB}\left (   \begin{array}{cccc}
  vw-c^{2}& \bb_{2}\ga_{2}-\al_{1}w&\al_{1}\ga_{1} - \bb_{2}v\\ 
\bb_{1}\ga_{1}-\al_{2}w& uw-b^{2}&\al_{2}\bb_{2}-  \ga_{1}u \\ 
\al_{2}\ga_{2}- \bb_{1}v&\al_{1}\bb_{1}-\ga_{2}u &uv-a^{2} \\ 
  \end{array} \right ) \label{1.2}
   \end{equation}
   Since $\BB^{-1}$ is an $M$ matrix   it's off diagonal terms are less than or equal to zero. This implies that
   \begin{equation}
 q_{1}:= { aw  \over bc}\ge1,\quad q_{2}:= { bv\over ac}\ge1,\quad  q_{3}:= {cu\over ab}\ge 1.\label{1.3}
   \end{equation}
We get this  by considering the pairs $\BB^{-1}_{1,2}$, $\BB^{-1}_{2,1}$ and  $\BB^{-1}_{1,3}$, $\BB^{-1}_{3,1}$ and  $\BB^{-1}_{2,3}$, $\BB^{-1}_{3,2}$ and using the hypothesis that $\BB^{-1}$ is an $M$ matrix.    
  We now note that 
  \bea
   {\det \AA\over abc}&=&\(\frac{uvw}{abc}-\frac{aw}{bc}\)+\( 1- {bv \over ac}\)+\(1-{ cu\over ab}\) \label{1.4a}\\\nn &=&\(q_{1}q_{2}q_{3}-q_{1}\)+\( 1- q_{2}\)+\(1-q_{3}\).
   \eea 
 One can check that on the range of $(q_{1},q_{2},q_{3})$ allowed in (\ref{1.3})
   $  \det \AA>0$  unless two of the $q_{i}=1,i=1,2,3$. (Set $q_{i}=1+\ep_{i},i=1,2,3$ and observe that   (\ref{1.4a}) is equal to $\ep_{1}\ep_{2}+\ep_{2}\ep_{3}+\ep_{1}\ep_{3}+\ep_{1}\ep_{2}\ep_{3}$.)  
   
Without loss of generality  suppose $q_{1}=q_{2}=1 $. 
In this case, when  $abc>0$,   $B^{-1}$ can not be an $M$-matrix.  We first note that if it is an $M$-matrix then $\{B^{-1}\}_{1,2}=\{B^{-1}\}_{2,1}=0$. This holds   since $\bb_{2}\ga_{2}\le\al_{1}w$ and $\bb_{1}\ga_{1}\le \al_{2}w$. Therefore, in order for $b^{2}c^{2}=a^{2}w^{2}$, which is the case when $q_{1}=1$, we must have $\bb_{2}\ga_{2}=\al_{1}w$ and $\bb_{1}\ga_{1}= \al_{2}w$.  A similar argument shows that  argument shows that $\{B^{-1}\}_{1,3}=\{B^{-1}\}_{3,1}=0$. In addition $ \{B^{-1}\}_{1,1}=0$, since $q_{1}=q_{2}=1 $ implies that $vw=c^{2}$. Consequently when $q_{1}=q_{2}=1 $ the first row of $B^{-1}=0$. But this is impossible since an $M$-matrix is invertible.

 We have now established that under the hypotheses of this theorem, when $abc>0$, $\det \AA>0$.
 
 Now suppose that $abc=0$. It is easy to see from (\ref{1.2}) that this implies that either two of these three numbers must be zero or all three must be zero. 
   For example, suppose that $\al_{1}=0$. Then by looking at $\BB^{-1}_{1,2}$ we see that either $\ga_{2}$ or $\bb_{2} $ are equal to zero.  Suppose $\al_{1}=\bb_{2} $=0. Then   
   \bea
    \det\BB 
  &=& u(vw-c^{2})+\al_{1}\bb_{1}\ga_{1}+\al_{2}\bb_{2}\ga_{2}-v\bb_{1}\bb_{2}-w\al_{1}\al_{2} \\
  &=& u(vw-c^{2}).\nn
   \eea
 In addition, since in this case $a=b=0$,   $\det \AA=u(vw-c^{2})$. Therefore, since $  \det\BB>0$ we have $  \det\AA>0$. A similar argument shows that when 
   $\al_{1}=\ga_{2} =0$ and $\BB^{-1}$ is an $M$-matrix, we also have $  \det\AA>0$.
     If all three numbers $a,b,c$ are equal to zero, $\det \AA=uvw$ and  so does  $\det\BB$, even though, once again, both $\BB$ and $\BB^{-1}$ can have non-zero off diagonal elements. In either case if $B^{-1}$ is an $M$-matrix,  $\det \AA>0$.
 
 Since 
 $\det\AA> 0$, 
we have
\begin{equation}
   \AA^{-1}=\frac{1}{\det \AA}\left (   \begin{array}{cccc}
  vw-c^{2}& bc-aw&ac-bv\\ 
bc-aw& uw-b^{2}& ab-cu \\ 
ac-bv& ab-cu & uv-a^{2} \\ 
  \end{array} \right ).
   \end{equation}
If $\AA^{-1}$ is not an $M$ matrix then one of the off diagonal terms of $ \AA^{-1}$ must be strictly positive. Without loss of generality we can suppose $aw<bc$. Then
\begin{equation}
   \al_{1}\al_{2}w^{2}<\bb_{1}\bb_{2}\ga_{1}\ga_{2}.
   \end{equation}
If this happens then one or both of the following inequalities must hold:
\begin{equation}
    \al_{1} w < \bb_{2} \ga_{2},\qquad  \al_{2}w <\bb_{1} \ga_{1}.
   \end{equation}
If either of these hold one of the off diagonal terms of $ \BB^{-1}$ must be positive. Therefore, $\BB^{-1}$ is not an $M$ matrix. This contradiction completes the proof. \qed

%Now suppose that one of the off diagonal elements of $\BB$ is zero. Then in order for $\BB^{-1}$ to be an $M$ matrix at least one other off diagonal element  of $\BB$ must be zero.  For example, if  $\al_{1}=0$ then at least one of $\bb_{2}$ or $\ga_{2}$ must be zero. Suppose $\bb_{2}=0$. Then
%\begin{equation}
%   \BB^{-1}=\frac{1}{det B}\left (   \begin{array}{cccc}
%  vw-c^{2}& 0&0\\ 
%\bb_{1}\ga_{1}-\al_{2}w& uw-b^{2}& -  \ga_{1}u \\ 
%\al_{2}\ga_{2}- \bb_{1}v& -\ga_{2}u &uv-a^{2} \\ 
%  \end{array} \right ) 
%   \end{equation}
%and
%\begin{equation}
%  \AA=\left (   \begin{array}{cccc}
%    u& 0&0\\ 
%   0& v& c  \\ 
%0&c & w \\ 
%  \end{array} \right )  . \end{equation}
%Thus if $\BB^{-1}$ is an $M$-matrix then $\AA^{-1}$ is an $M$-matrix. 

%It is also possible that $\BB^{-1}$ is an $M$-matrix and  $\AA$ is a diagonal matrix. In this case, once again, $\AA^{-1}$ is an $M$-matrix. \qed

\bl \label{cor-1.1} Let $\{\th_{x},x\in R^{1}\}$ be a permanental process with kernel $\Ga(x,y)$. Let
\begin{equation}
   d(x,y):=\(\Ga(x,x)+\Ga(y,y)-2(\Ga(x,y)\Ga(y,x))^{1/2}\)^{1/2}.
   \end{equation}
Suppose that for any three different real numbers, $x_{1}$, $x_{2}$, $x_{3}$,  the $3\times 3$ matrix $\wt\Ga:=\{\Ga(x_{i},x_{j})\}_{i,j=1,2,3}$ has an inverse that is an $M$-matrix. Then $d(x,y)$ is a metric on $R^{1}$.
 \el
 
 \Proof For any $x_{1}$, $x_{2}$, $x_{3}$ in $R^{1}$ it is easy to see that  
 \begin{equation}
    d(x_{1},x_{2}) \le d(x_{1},x_{3})+ d(x_{3},x_{2}).\label{1.9}
   \end{equation}
   This follows because by  Lemma \ref{lem-1.1}   the determinant of  
 \begin{equation}
\!\!\!\left (   \!\!\begin{array}{cccc}
    \Ga(x_1,x_1)& (\Ga(x_1,x_2)\Ga(x_2,x_1))^{1/2}& (\Ga(x_1,x_3)\Ga(x_3,x_1))^{1/2}\\ 
   (\Ga(x_1,x_2)\Ga(x_2,x_1))^{1/2}& \Ga(x_2,x_2)& (\Ga(x_{3},x_2)\Ga(x_2,x_{3}))^{1/2}  \\ 
(\Ga(x_1,x_3)\Ga(x_3,x_1))^{1/2}&(\Ga(x_{3},x_2)\Ga(x_2,x_{3}))^{1/2} & \Ga(x_3,x_3) \\ 
  \end{array}\!\! \right )
   \end{equation}
 is strictly positive.  The determinant of the   principle 2$\times$2 minor is also strictly positive by \cite[(1.5)]{perm}. Therefore the matrix is positive definite and thus it is the covariance matrix of a Gaussian random variable which we denote by $(G_{x_{1}},G_{x_{2}},G_{x_{3}})$. Since   \begin{equation}
       d(x_{i},x_{j})=\(E(G_{x_{i}}-G_{x_{j}})^{2}\)^{1/2}\qquad i,j=1,2,3,
   \end{equation}
 we get (\ref{1.9}).\qed
 
 \begin{remark} {\rm Note that $\det \BB-\det \AA=(\al_{1}\bb_{1}\ga_{1}-\al_{2}\bb_{2}\ga_{2})^{2}$ so it is possible to have $\det\BB>0$ and $\det \AA\le 0$, but not when $\BB^{-1}$ is an $M$-matrix.

 }\end{remark}
 
 \subsection{Continuity of associated permanental processes}
 
\medskip	 By Theorem \ref{theo-4.1} when  $u(x,y)$, $x,y\in T$, is the potential density of a transient Markov process with state space $T$, a $\bb$-permanental process $\th_{\bb}= \{\th_{x,\bb}: x\in T\}$, with kernel $u$, exists for all $\bb>0$. In addition, in this case, by  Lemma \ref{lem-4.2}, the function $d$ in (\ref{1.2w}) is a metric on $T$. Using these facts we get the following extension of Theorem \ref{theo-1.1a}:
  
   \bt\label{theo-4.2}Ê Ê  Let $\th_{\bb}=\{\th_{x,\bb}, x\in T\}$ be a $\bb$-permanental process associated with a transient Markov process with potential $u(x,y)$ satisfying $\sup_{x\in T} u(x,x)<\ff $.  Let $\wt d$ be as given in (\ref{4.2w}) and let $D$ denote the $\wt d$ diameter of $T$.  Assume that   there exists
a probability measure $\mu$ÊÊ Ê on $\mathcal{B}(T, \wt d)$ such that  
\be
J_{ \wt  d }(D)<\ff\label{1.1v}.
\ee 
Then there exists a version $\th'_{\bb}=\{\th'_{x,\bb},x\in T\}$ of $\th_{\bb}$ such that for any $x_{0}\in T$
 \be
 \|\sup_{x\in T}\th' _{x}  \|_{\psi_{1}}\le   C_{\bb}\(\| \th' _{x_{0}}  \|_{\psi_{1}}\label{1.5mm} +   \sup_{x\in T}\Ga(x,x)  
J^{2}_{\wt  d}(D)\),   
 \ee where $C_{\bb}$ is a constant.   
If
\begin{equation}
   \lim_{\de\to 0}J_{ \wt d}(\de)=0,\label{1.8h}
   \end{equation}
there exists a version  $\th'_{\bb}=\{\th'_{x,\bb},x\in T\}$ of $\th_{\bb}$  that is continuous on $(T,\wt d)$. 

If (\ref{1.8h}) holds and 
\begin{equation} 
 \lim_{\de\to 0}{J_{\wt  d}(\de)\over \de}=\ff
  \label{2.18v},
   \end{equation}
   then      
   \begin{equation}
   \lim_{\de\to0}   \sup_{\stackrel{x,y\in T}{  \wt d (x,y)\le \de}}\frac{|\th' _{  x,\bb} -\th' _{y,\bb}|}{J_{  \wt d} (  \wt d(x,y) /2) }\le C_{\bb}\(\sup_{x\in T}\th'_{x,\bb}\)^{1/2}\quad a.s.,\label{2.1svqq}
   \end{equation} where $C_{\bb}$ is a constant depending only on $\bb$.
   \et

\Proof  Given $\bb$ choose $p$ so that $\bb\le p/2$. We have 
 \be 
 \({\th_{x,\bb} } -{\th_{y,\bb}  }\)+ \({\th'_{x,p/2-\bb} } -{\th'_{y,p/2-\bb}  }\) \stl  \sum_{j=1}^{p}\({\th_{x,1/2} }^{(j)}-{\th_{y,1/2} }^{(j)} \), \label{4.20}  \ee  
where $\th_{\bb}$ and $\th_{\bb'}$ are independent and   where the notation $(j)$ denotes independent copies. Let $\|\cd\|$ be a pseudo norm on real valued functions on $T$. By (\ref{4.20}),  
 \be E\| \({\th_{x,\bb} } -{\th_{y,\bb}  }\)+ \({\th'_{x,p/2-\bb} } -{\th'_{y,p/2-\bb}  }\) \|\\
 \le p E \|\th_{x,1/2} -\th_{y,1/2} \|.     \ee 
 Using the same argument we used in (\ref{3.26})--(\ref{3.29}) we get
 \be  
   E\|{\th_{x,\bb} } -{\th_{y,\bb} } \|\le 2(p-\bb) E\| {\th_{x,1/2} }  -{\th_{y,1/2} } \|.\label{3.35}
 \ee 

  Similarly, we also note that starting with
  \be 
 {\th_{x,\bb} } +  {\th'_{x,p/2-\bb} }   \stl  \sum_{j=1}^{p}{\th_{x,1/2} }^{(j)}  \label{4.20a}  \ee  
 we can show that
 \be  
   E\|{\th_{x,\bb} } \|_{\psi_{1}}\le 2(p-\bb) E\| {\th_{x,1/2} }    \|_{\psi_{1}}.\label{3.35a}
 \ee
Therefore, (\ref{1.5mm}) follows from (\ref{1.5m}).

Using (\ref{3.35}) we have
 \be  
   E \sup_{\stackrel{s,t\in T'}{  \wt d(s,t)\le \de}}|{\th_{x,\bb} } -{\th_{y,\bb} } |\le 2(p-\bb) E\sup_{\stackrel{s,t\in T}{   \wt d(s,t)\le \de}}| {\th_{x,1/2} }  -{\th_{y,1/2} } |\label{3.35qq}
 \ee 
 for all finite subsets $T'$ of $T$. By (\ref{1.8h}) and (\ref{3.6.90hhv})   there exists a version $\th'_{\bb}=\{\th'_{x,\bb},x\in T\}$ of $\th_{\bb}$ \be  
 \lim_{\de\to 0}  E \sup_{\stackrel{s,t\in T'}{   \wt d(s,t)\le \de}}|{\th'_{x,\bb} } -{\th'_{y,\bb} } |=0 \ee
 This implies that $\th_{\bb}'$    is continuous on $(T,\wt d)$.
 
 By (\ref{2.18v}) and (\ref{2.1sv})  there exists a positive function $g(\de)$, satisfying $\lim_{\de\to 0}g(\de)=0$, such that
   \begin{equation}
 P\(   \sup_{\stackrel{s,t\in T}{\wt  d (s,t)\le \de }}\frac{|\th' _{  s,1/2} -\th' _{t,1/2}|}{J_{ \wt d} (  \de  /2) }>31\(\sup_{x\in T}\th'_{x,1/2}\)^{1/2}\)<g(\de) .\label{2.1svqq1}
   \end{equation}
 Using the argument in (\ref{3.37}) we get that  there exists a   constant $C_{\bb}$ and a positive function $g'(\de)$, satisfying $\lim_{\de\to 0}g'(\de)=0$, such that
 \begin{equation}
 P\(   \sup_{\stackrel{s,t\in T}{ \wt d (s,t)\le \de }}\frac{|\th' _{  s,\bb} -\th' _{t,\bb}|}{J_{ \wt d} (  \de  /2) }>C_{\bb}\(\sup_{x\in T}\th'_{x,1/2}\)^{1/2}\)<g'(\de).\label{2.1svqq2}
   \end{equation}
 This implies that \begin{equation}
   \lim_{\de\to0}   \sup_{\stackrel{x,y\in T}{  \wt d (x,y)\le \de}}\frac{|\th' _{  x,\bb} -\th' _{y,\bb}|}{J_{  \wt d} (  \de/2) }\le C_{\bb}\(\sup_{x\in T}\th'_{x,\bb}\)^{1/2}\quad a.s.\label{2.1svqq3}\ee
 By \cite[Lemma 7.1.6]{book} this is equivalent to (\ref{2.1svqq}).\qed

      \section{Permanental processes associated with L\'evy processes}\label{sec-pfbm}  

  It is interesting to see what kernels    of permanental processes, that are not Gaussian squares, look like. We can do this by finding 
  the kernels of    permanental processes associated with  L\'evy processes. Let $X=\{X_{ t},t\in R_{+}\}$ be a L\'evy process with characteristic function
\begin{equation}
 Ee^{i\la X}=e^{-\psi(\la)t}.\label{10.1}
 \end{equation}
 We refer to $\psi$ as the characteristic exponent of $X$.

We assume that $X$ has local times $\{L_{t}^{x}, (x,t)\in R\times R_{+}\}$. Set 
\begin{equation}
 u_{T_{0}}(x,y)=E^{x}\(L_{T_{0}}^{y}\),\label{10.2q}
 \end{equation}
 where $T_{0}$ is the first hitting time of $X$ at zero. We consider the transient Markov process $\wt X=\{\wt X_{t} \}$   that is $X$ killed at the first time  it hits zero. The function $u_{T_{0}}(x,y)$ is the zero potential of $\wt X$ and thus is also the kernel of a permanental process.

\begin{lemma} \label{lem-10.1}
\begin{equation}
 u_{T_{0}}(x,y)=R(x,y)+H(x,y) \quad\mbox{and}\quad  u_{T_{0}}(y,x)=R(x,y)-H(x,y) \label{10.2}
 \end{equation}
 where
\begin{equation}
 R(x,y)=    \frac1{ \pi} \int_{0}^\ff\frac{\(1- \cos\la x -\cos\la y+\cos\la(x-y) \)\RR e\,\psi( \la)}{|\psi(\la)|^{2}}\,d\la\label{10.3}
 \end{equation}
 and
\begin{equation}
H(x,y)= -\frac1{ \pi} \int_{0}^\ff\frac{ (\sin\la x-\sin\la y-\sin\la(x-y))\,  \II m\,\psi( \la)}{|\psi(\la)|^{2}}\,d\la.\label{10.4}
 \end{equation}  
\end{lemma}

\Proof   It follows from \cite[Theorem 4.2.4]{book}, modified for  non-symmetric L\'{e}vy processes, that for $x,y\in R^1$
   \be
         u_{T_0}(x,y)=  \phi(x)+\phi(-y)-\phi(x-y),  \label{mp11.6}
        \ee
        where
         \bea
        \phi(x )&=&\frac1{2\pi}\int_{-\ff}^\ff\frac{1- e^{i \la x }}{\psi(\la)}\,d\la\label{mp11.7}\\
        &=&\nn\frac1{2\pi}\int_{-\ff}^\ff\frac{\(1- e^{i \la x }\)\psi(-\la)}{|\psi(\la)|^{2}}\,d\la \\
         &=&\nn\frac1{2\pi}\int_{0}^\ff\frac{\(1- e^{i \la x }\)\psi(-\la)}{|\psi(\la)|^{2}}\,d\la +\frac1{2\pi}\int_{0}^\ff\frac{\(1- e^{-i \la x }\)\psi( \la)}{|\psi(\la)|^{2}}\,d\la\\
     &=&    \frac1{ \pi}\RR e\(\int_{0}^\ff\frac{\(1- e^{i \la x }\)\psi(-\la)}{|\psi(\la)|^{2}}\,d\la \)\nn
     \eea
Since $\ov{\psi(\la)}={\psi(-\la)}$ we see that  
     \be\phi(x)=   \frac1{ \pi}\(\int_{0}^\ff\frac{\(1- \cos\la x  \)\RR e\,\psi( \la)}{|\psi(\la)|^{2}}\,d\la-\int_{0}^\ff\frac{  \sin\la x\,  \II m\,\psi( \la)}{|\psi(\la)|^{2}}\,d\la \).      \label{10.9}     \ee 
Therefore,          
\be   { \phi(x) +\phi(-x)\over 2}=\frac1{ \pi} \int_{0}^\ff\frac{\(1- \cos\la x  \)\RR e\,\psi( \la)}{|\psi(\la)|^{2}}\,d\la. \label{10.10qq} \ee    
and
\be   { \phi(x) -\phi(-x)\over 2}=-\frac1{ \pi} \int_{0}^\ff\frac{ \sin\la x\,  \II m\,\psi( \la)}{|\psi(\la)|^{2}}\,d\la. \label{10.11} \ee    

\medskip	We write 
\begin{equation}
 u_{T_{0}}(x,y)=\frac{u_{T_{0}}(x,y)+u_{T_{0}}( y,x)}{2}+\frac{u_{T_{0}}(x,y)-u_{T_{0}}( y,x)}{2}\label{10.6}
 \end{equation}
   and
 \begin{equation}
 u_{T_{0}}(y,x)=\frac{u_{T_{0}}(x,y)+u_{T_{0}}( y,x)}{2}-\frac{u_{T_{0}}(x,y)-u_{T_{0}}( y,x)}{2}.\label{10.7}
 \end{equation}  
 Set
 \begin{equation}
 R(x,y)=\frac{u_{T_{0}}(x,y)+u_{T_{0}}( y,x)}{2}\quad \mbox{and}\quad H(x,y)=\frac{u_{T_{0}}(x,y)-u_{T_{0}}( y,x)}{2}.
 \end{equation}
Using (\ref{mp11.6}) and the following equations we get (\ref{10.2})--(\ref{10.4}).\qed

We now consider the transient Markov process $\ov X=\{\ov X_{t} \}$    that is $X$ killed at   $\xi_{\al}$, an independent exponential time with mean $\al>0 $ . We assume that $X$ has local times $\{L_{t}^{x}, (x,t)\in R\times R_{+}\}$ and set 
\begin{equation}
 u^{\al} (x,y)=E^{x}\(L_{\xi_{\al}}^{y}\) \label{10.2qq}.
 \end{equation}
The function $u^{ \al}(x,y)$ is the zero potential of $\ov X$ and thus is also the kernel of a permanental process. 

\begin{lemma} \label{lem-10.2} For $\al>0$, the $\al$-potential density
\begin{equation}
 u^{ \al}(x,y)=  R_{\al}(x,y)+  H_{\al}(x,y) \quad\mbox{and}\quad  u ^{\al}(y,x)=R_{\al}(x,y)-  H_{\al}(x,y)\label{10.2ss},
 \end{equation}
 where
\begin{equation}
 R_{\al}(x,y)=    \frac1{ \pi} \int_{0}^\ff\frac{ \cos\la(x-y)  \RR e(\al+\psi( \la))}{|\al+\psi(\la)|^{2}}\,d\la\label{10.3ss}
 \end{equation}
 and
\begin{equation}
 H_{\al}(x,y)=  \frac1{ \pi} \int_{0}^\ff\frac{  \sin\la(x-y) \,  \II m(\al+\psi( \la))}{|\al+\psi(\la)|^{2}}\,d\la.\label{10.4ss}
 \end{equation}  
\end{lemma}

\Proof  By (\ref{10.1}) 
\begin{equation}
  \int_{-\ff}^{\ff }e^{i\la y}p_{t}(0,y)\,dy=e^{-t\psi(\la)}
  \end{equation}
Therefore
\begin{equation}
  p_{t}(0,y)=\frac{1}{2\pi}   \int_{-\ff}^{\ff }e^{-i\la y} e^{-t\psi(\la)}\,d\la
  \end{equation}
and, similarly to (\ref{mp11.7})
 \bea
 u^{\al}(0,y)&=&\int_{0}^{\ff}e^{-\al t} p_{t}(0,y)\,dt=\frac{1}{2\pi}   \int_{-\ff}^{\ff }{e^{-i\la y} \over \al+\psi(\la)} \,d\la\\ \nn&=&   \frac1{ \pi}\RR e\(\int_{0}^\ff\frac{  e^{ i \la y } (\al+\psi( \la))} {|\al+ \psi(\la)|^{2}}\,d\la \) \\
&= &  \frac1{ \pi}\(\int_{0}^\ff \frac{ \cos\la y \RR e\,(\al+\psi( \la))}{|\al+ \psi(\la)|^{2}}\,d\la\right.\nn\\
&&\hspace{1in}\left.-\int_{0}^\ff\frac{  \sin\la y\,  \II m\,(\al+\psi( \la))}{|\al+ \psi(\la)|^{2}}\,d\la \).    \nn        \eea 
Using the fact that $ u^{\al}(x,y) = u^{\al}(0,y-x)$      we see that   
\be   {  u^{\al}(x,y)+ u^{\al}(y,x)\over 2}=\frac1{ \pi} \int_{0}^\ff\frac{  \cos\la (x-y)  \RR e\,(\al+\psi( \la))}{|\al+ \psi(\la)|^{2}}\,d\la. \label{10.10} \ee    
and
\be  {  u^{\al}(x,y)- u^{\al}(y,x)\over 2}= \frac1{ \pi} \int_{0}^\ff\frac{ \sin\la( x-y)\,  \II m\,(\al+\psi( \la)}{|\al+ \psi(\la)|^{2}}\,d\la,\label{10.11s} \ee    
which gives (\ref{10.2ss}).\qed

\begin{remark} \label{rem-9.1}{\rm Denote  $R$ and $H$ in (\ref{10.3}) and (\ref{10.4}) by $R_{0}$ and $H_{0}$. Since 
 \begin{equation}
 1- \cos\la x -\cos\la y+\cos\la(x-y)=( 1- \cos\la x )( 1- \cos\la y )+\sin\la x\sin\la y
 \end{equation}
 and 
   \begin{equation}
\cos\la(x-y)=  \cos\la x   \cos\la y  +\sin\la x\sin\la y
 \end{equation}
we see that $R_{\al}(x,y)$ is positive definite for all $\al\ge 0$. Therefore they are  the covariances of  mean zero Gaussian processes  $G_{\al}=\{G_{\al}(x),x\in R^{1}\}$. It is easy to see that for $\al>0$, $G_{\al}$ is a stationary Gaussian process and that $G_{0}$, has stationary increments, and that $G_{0}(0)=0$. Furthermore    
\begin{equation}
 E(G_{0}(x)-G_{0}(y))^{2}= \frac2{ \pi} \int_{0}^\ff\frac{\(1- \cos\la( x-y)   \)\RR e\,\psi( \la)}{|\psi(\la)|^{2}}\,d\la 
 \end{equation}
 and when  $E(G^{2}(x))=1$,
  \begin{equation}
 E(G_{\al}(x)-G_{\al}(y))^{2}= \frac2{ \pi} \int_{0}^\ff\frac{ (1-  \cos\la( x-y)   )\RR e(\al+\,\psi( \la))}{|\al+ \psi(\la)|^{2}}\,d\la ,\quad\forall \al>0.
 \end{equation} 
Since  $R_{\al}(x,y)$ is positive definite for all $\al\ge 0$, it follows from from   (\ref{10.2}) that $u_{T_{0}}(x,y)$ and $u^{\al}(x,y)$, $\al>0$,  are positive definite. This is easy to see since
\begin{equation}
  u_{T_{0}}(x,y)=R_{0}(x,y)+H_{0}(x,y)
  \end{equation}
  and $H_{0}(x,x)=0$, and $H_{0}(x,y)=-H_{0}(y,x)$, and similarly for $u^{\al}(x,y)$.

   \medskip	 
   Set \begin{equation}
 d_{G_{\al}}(x,y)= \(E(G_{\al}(x)-G_{\al}(y))^{2}\)^{1/2}\quad\al\ge0.
 \end{equation}
  We show in \cite[Lemma 5.5]{perm} that this metric is equivalent to   $\wt d$ in (\ref{4.2w}). Therefore it follows from   Theorem \ref{theo-4.2}, and \cite[Chapter 6]{book}, that the   denominator in (\ref{2.1svqq})  for a   permanental processes with kernel $u_{T_{0}}(x,y)$ or $u^{\al}(x,y)$, $\al>0$, is the same, up to a constant multiple,  as the uniform  modulus of continuity of Gaussian processes with increments variance $d_{G_{\al}}$.
}\end{remark}

    \subsection{$\bf FBMQ^{\al ,\bb}$}   \label{subsec-5.1}

 There is a great deal of interest in fractional Brownian motion. 
Therefore, there should be considerable interest in permanental process for which the function $R(x,y)$ in  (\ref{10.3}) is the covariance of fractional Brownian motion. To study these processes we use Lemma \ref{lem-10.1} to find the kernel of the permanental process associated with the non-symmetric stable process of index   $1< \al+1  <2$. This process has a L\'{e}vy exponent of the form
\begin{equation}
\psi (\la) =c|\la|^{\al+1}(1-i\bb \mbox{ sign }(\la)\,\tan ((\al+1)\pi/2))\label{9.1}
\end{equation}
where $c$ is an arbitrary constant which we set to 1, and $|\bb|\leq 1$.  We refer to this permanental process   as $\mbox{FBMQ}^{\al ,\bb}$.
We use this notation because when $\bb=0$, $\psi(\la)$ is the L\'evy exponent of   the  L\'evy  process associated with  fractional Brownian motion of index $\al$, (i.e. FBM). We add the $Q$, for quadratic, to denote the square of this process, as one does in the designation of the squared Bessel processes, (BESQ).  

 \begin{lemma} Let   $X=\{X_{ t},t\in R_{+}\}$ be a L\'evy process with characteristic given by (\ref{9.1}). Let $u_{T_{0}}(x,y)$ be as given in (\ref{10.2q}) and write $u_{T_{0}}(x,y)$ and $u_{T_{0}}(y,x)$ as in (\ref{10.2}).  In this case
 \begin{equation}
 R(x,y)=
C_{\al,\bb}  \(|x|^{ \al}+|y|^{ \al}-
|x-y|^{ \al}\),     \label{9.6}
\ee
and
       \begin{eqnarray}
\lefteqn{H(x,y) }\label{9.6a} \\
&& =\bb C_{\al,\bb}  \(\mbox{ sign }(x)|x|^{ \al}-\mbox{ sign }( y)|y|^{ \al}-
\mbox{ sign }(x-y)|x-y|^{ \al}\),\nn
\end{eqnarray}
where 
\begin{equation}
 C_{\al,\bb}={-\sin \left(( \al+1)\frac{\pi}{2} \right) 
\Gamma( -\al ) \over \pi (1+\bb^{2}  \,\tan^{2} ((\al+1)\pi/2))}>0.
 \end{equation}
 Consequently
\begin{eqnarray}
\lefteqn{ u_{T_{0}}(x,y)=
C_{\al,\bb}  \(\(1-\bb \mbox{ sign }(x)\)|x|^{ \al}\right.
\label{sas}}\\
&&\qquad\quad \left. +\(1-\bb \mbox{ sign }(-y)\)|y|^{ \al}-
\(1-\bb \mbox{ sign }(x-y)\)|x-y|^{ \al}\)  \nonumber.
\end{eqnarray}
\end{lemma}

\Proof  Let  $c_{\al,\bb}= (1+\bb^{2}  \,\tan^{2} ((\al+1)\pi/2))^{-1}$. Then
\be 
{ \(1- \cos\la x  \)\RR e\,\psi( \la)\over |\psi(\la)|^{2}}={c_{\al,\bb}\(1- \cos\la x  \)\over |\la|^{ \al+1 }}
\ee
and
\be  \frac{  \sin\la x\,  \II m\,\psi( \la)}{|\psi(\la)|^{2}}=-{ c_{\al,\bb}\bb \mbox{ sign }(\la)\,\tan ((\al+1)\pi/2)\sin\la x\over |\la|^{\al+1}}
\end{equation}
Therefore, by a change of variables
\bea
&& \int_{0}^\ff\frac{  \sin\la x\,  \II m\,\psi( \la)}{|\psi(\la)|^{2}}\,d\la\label{10.21}\\
 & &\qquad=-c_{\al,\bb}\,\bb\tan ((\al+1)\pi/2)    \int_{0}^\ff{\sin\la x \over   \la ^{\al+1}}\,d\la\nn\\
 & &\qquad=-c_{\al,\bb}\,\bb\tan ((\al+1)\pi/2)  |x|^{\al }\mbox{   sign }( x )   \int_{0}^\ff{\sin\la \over    \la ^{\al+1}}\,d\la\nn 
 \eea
and
\be  \int_{0}^\ff { \(1- \cos\la x  \)\RR e\,\psi( \la)\over |\psi(\la)|^{2}}\,d\la\\
= c_{\al,\bb}  |x|^{\al  } \int_{0}^\ff{\(1- \cos\la    \) \over    \la ^{\al+1}}\,d\la  .
 \ee 
It follows from \cite[17.33.2]{GR} and an integration  that
\begin{equation}
\int_{0}^\ff\frac{ 1- e^{i \la} }
       { \la ^{\al+1}}\,d\la =-\( \sin \left( (\al+1)\frac{\pi}{2} \right)+i \cos \left( (\al+1)\frac{\pi}{2} \right)\)
\Gamma( -\al ).\label{9.5}
\end{equation} 
Combining  (\ref{10.21})--(\ref{9.5}) we see that
\begin{equation}
\frac{1}{\pi} \int_{0}^\ff { \(1- \cos\la x  \)\RR e\,\psi( \la)\over |\psi(\la)|^{2}}\,d\la=   C_{\al,\bb}\, |x|^{\al } \label{9.8}
\end{equation}
and 
\begin{equation}
\frac{1}{\pi}\int_{0}^\ff\frac{  \sin\la x\,  \II m\,\psi( \la)}{|\psi(\la)|^{2}}\,d\la=   \bb\mbox{ sign }(x)  C_{\al,\bb}\, |x|^{\al }.\label{9.8a}
\end{equation} 
Using (\ref{9.8}) and  (\ref{9.8a}) in (\ref{10.3}) and (\ref{10.4}) we get  (\ref{9.6}) and (\ref{9.6a}).\qed

\begin{remark} {\rm It follows from Remark \ref{rem-9.1} and \cite[(7.186)]{book} that when $\th_{x}$ is $\mbox{FBMQ}^{\al ,\bb}$ then for any $T<\ff$
     \begin{equation}
 \lim_{\de\to0}   \sup_{\stackrel{x,y\in [0,T]}{  d (x,y)\le \de}}\frac{|\th_{  x} -\th _{y}|}{( |x-y|^{\al }\log 1/ |x-y| )^{1/2} }\le C_{T}\(\sup_{x\in [0,T]}\th_x\)^{1/2}\quad a.s.\label{2.1fbm}
 \end{equation}
In addition by \cite[Theorem 4.2]{perm} and \cite[Example 7.6.6, (1)]{book} for $x_{0}\ne 0$
\begin{equation}
 \lim_{\de\to0}   \sup_{ \hat d(x,x_{0})< \de/2 }\frac{ |\th _{x}-\th _{x_{0}}|}{( |x-x_{0}|^{\al} \log \log1/ |x-x_{0}| )^{1/2}}\le C\, \th^{1/2}_{x_{0}}  \qquad a.s.\label{2.1wj}
 \end{equation}
 and
 \begin{equation}
   \lim_{\de\to0}   \sup_{ \hat d(x,0)< \de/2 }\frac{ \th _{x} }{  |x |^{\al } \log \log1/ |x |  }\le C   \qquad a.s. 
  \end{equation}
}
\end{remark}
 
A very interesting unanswered question is: can we find $\bf FBMQ^{\al }$ for $1\le\al<2$ other than the square of a fractional Brownian motion itself? We give some examples to explain what we have in mind. Let $\wt G_{\al}=\{\wt G_{\al}(x),x\in R_{+}\}$ be fractional Brownian motion  with index $1\le\al<2$, i.e. $\wt G_{\al}(0)=0$ and 
\begin{equation}
\si^{2}(x-y):=    E(\wt G_{\al}(x)-\wt G_{\al}(y))^{2}=|x-y|^{\al}.\label{5.44}
   \end{equation} 
 It follows that the covariance
 \begin{equation}
    \wt \Ga(x,y):=    E(\wt G_{\al}(x) \wt G_{\al}(y))=\frac{1}{2}\(|x|^{\al}+|y|^{\al}-(|x-y|^{\al})\).\label{6.17}
   \end{equation}
   Consider the kernel
   \begin{equation}
   \GG(x,y)=\left \{
  \begin{array}{  ll  }
      \wt \Ga(x,y)+\ep\,\si (x-y) &x\ge y \\ \\
      \wt \Ga(x,y)-\ep\,\si (x-y) &x<y 
  \end{array}\right.\label{5.47}
    \end{equation}
where $\ep>0$ is some small number. To simplify the notation set $\ov d(x,y)=d(x,y)/4\sqrt{2/3}$. We can write
\bea
  \ov d^{2}(x,y)&=&   \GG(x,x)  +\GG(y,y) - \GG(x,y) - \GG(y,x) \label{6.19}\\\nn
&&\qquad+    \GG(x,y) + \GG(y,x)-2\(\GG(x,y)  \GG(y,x)\)^{1/2} \\
&=&   \wt \Ga(x,x)  + \wt \Ga(y,y)  -  2\wt \Ga(x,y) +   ( \GG^{1/2}(x,y) -\GG^{1/2}(y,x))^{2} \nn \\
&\le& \si^{2}(x-y) +   | \GG (x,y) -\GG (y,x)|\nn\\
&\le&(1+2\ep) \si^{2}(x-y) .\nn
     \eea
Also we see from the second equality in (\ref{6.19}) that $  \ov d^{2}(x,y)\ge \si^{2}(x-y) $. 

We do not know whether $\GG(x,y)$ is the kernel of a permanental process. If it is, it would be a generalization of the square of fractional Brownian motion. It would be very useful to know whether such generalizations exist. 

Here is a very simple example of the above. Consider (\ref{6.17}) when $\al=1$. This is the covariance of Brownian motion. Consider $\GG(x,y)$  in this case for three values, $0<x<y<z$. The kernel obtained is 
\be \AA_{\ep}:=   \left (   \begin{array}{cccc}
    x& x-\ep|x-y|^{1/2}&  x-\ep|x-y|^{1/2}\\ 
   x+\ep|x-y|^{1/2}& y&  y-\ep|x-y|^{1/2}  \\ 
  x+\ep|x-y|^{1/2}& y+\ep|x-y|^{1/2} & z \\ 
  \end{array} \right ) .
  \ee
When $\ep=0$, $\AA_{0}$ is the kernel of $\{B^{2}_{x},B_{y}^{2}, B_{z}^{2}\}$ where 	$\{B _{x},B_{y} , B_{z} \}$ are three values of standard Brownian motion. One can check that $\AA_{0}^{-1}$ in an $M$-matrix. This is not a surprise since we now that Brownian motion is an associated process. (See e.g. \cite[Lemma 2.5.1]{book}.) However,  $\AA_{\ep}^{-1}$ is not an $M$-matrix for any $\ep\ne 0$. It would be interesting to know whether $\AA_{\ep}$ can be the kernel of a permanental vector.

\begin{remark} {\rm One may consider kernels of the form (\ref{5.47}) for any Gaussian process with stationary increments and try to answer the question, ``Are they kernels of permanental processes?''

 }\end{remark}

\section{Some  problems} \label{sec-6}

\noindent {\bf 1. Necessary and sufficient conditions for the continuity of permanental processes.}   This question could be phrased more specifically as, ``Is Theorem \ref{theo-1.1a} a necessary and sufficient condition  for the continuity of 1/2-permanental processes?'' It is clear that Theorem \ref{theo-1.1a} is best possible because it is a necessary and sufficient condition  for the continuity of the square of a Gaussian process with continuous covariance; (see \cite[Section 7.2]{perm}).  

When one considers this question one is immediately aware of  how many of the properties of  Gaussian processes that are used in the proof of a necessary condition for continuity or boundedness are unknown for permanental processes. Let $\th=\{\th_{x}, x\in T\}$ be a $1/2$-permanental process. If $\sup_{x\in T} \th_{x} <\ff$, is $E(\sup_{x\in T} \th_{x} <\ff)$?  What about $\|\sup_{x\in T} \th_{x}  \|_{\psi_{1}}$? If the metric $d$ in (\ref{1.2w}) of one permanental process is uniformly less that the metric $d$ of another permanental process, does the continuity or boundedness of the process with the larger metric imply the continuity or boundedness of the process with the smaller metric? More generally, are there any analogues of Slepian's Lemma for permanental processes? Gaussian processes often satisfy many important  zero one laws.
What zero one laws do permanental processes satisfy? All of these are interesting questions.

\medskip	\noindent {\bf 2. Independence of  the components of a permanental vector.} It follows from Lemma \ref{lem-1} that if 
  $\th:=\{\th_{x}, x\in T\} $ is a  1/2-permanental process
with kernel $\Ga$ then for any pair $x , y$, 
\begin{equation}
\mbox{cov}   (\th_{x},\th_{y})=2\Ga(x,y)\Ga(y,x) \label{funda}
   \end{equation}
If $\Ga(x,y)\Ga(y,x)=0$ it  follows from Lemma \ref{lem-1} that $\th_{x}$ and $\th_{y}$ are independent. Another interesting question is, ``If the elements of   $ \{\th_{x}, x\in T\} $ are pairwise independent are they independent?''

Suppose that $T=\{1,2,3\}$ and the kernel of $\th$ is 
\be  \(   \begin{array}{cccc}
    1& a& 0 \\ 
   0& 1& b  \\ 
 c& 0 & 1    
  \end{array}  \).\label{6.2}
\ee
where $abc\ne 0$, and otherwise $a,b,c$ are arbitrary. 

By the general formula for the moments of $\bb$-permanental processes
\bea
E\( \prod_{j=1}^{n}\th_{x_{j}}  \)&=& \mbox{Perm}_{\bb}\(\lc \Ga (x_{i},x_{j})\rc_{1\leq i,j\leq n}\)\label{fs.14}\\
&=& \sum_{\pi}\al^{c(\pi)} \prod_{j=1}^{n}\Ga (x_{j},x_{\pi(j)})\nn
\eea
where $c(\pi)$ is the number of cycles in the permutation $\pi$ of $[1,n]$.
  (See \cite[Proposition 4.2]{VJ}.) Therefore,    if the   kernel given by 
(\ref{6.2}) defines a $1/2$-permanental   vector $(\th_{1},\th_{2},\th_{3})$, then we would have $E(\th_{1}\th_{2}\th_{3})$ =1/8+(1/2)abc. Whereas, if $ \th_{1},\th_{2}, \th_{3} $ were independent   $E(\th_{1}\newline\th_{2}\th_{3})=1/8$. On the other hand by (\ref{funda})  $ \th_{1},\th_{2},\th_{3} $ are pairwise independent.

The kernel in (\ref{6.2}) and its transpose are the only possibilities for   $\th_{1},\th_{2} ,\th_{3}$ to be pairwise independent but not independent. So in this case the problem becomes, can this be the kernel of a   $1/2$-permanental vector. We can show that it can not. But in $R_{+}^{4}$ we can not answer this question. In   $R_{+}^{4}$ we can have kernels of the form
\be  \(   \begin{array}{cccc}
    1& a& 0&0 \\ 
   0& 1& b &e \\ 
 0& 0 & 1 &f \\ 
d &0&0&1
  \end{array}  \).\label{6.2qq{perm}}
\ee

 \medskip	\noindent {\bf 3. Permanental vectors in $R_{+}^{3}$.} Given a $3\times 3$ matrix with components that satisfy (\ref{vj.1})--(\ref{vj.2}), is it the kernel of a permanental process? If the inverse of the matrix is an $M$-matrix, the answer is yes. At this point we don't know of any examples of matrices that are kernels of a permanental process with  inverses that are   not  $M$-matrices, except, of course, when the permanental process is the square of a Gaussian process. We think that there are many examples. The class of associated vectors that are Gaussian squares is a very small subset of all vectors that are Gaussian squares. Roughly speaking, they don't include any vectors that are squares of Gaussian processes that are smoother than Brownian motion.  If, again, roughly speaking, one can always find permanental vectors with kernels that deviate slightly from the kernels of Gaussian squares that are infinitely divisible, one should also be able to find permanental vectors with kernels that deviate slightly from the kernels of Gaussian squares that are not infinitely divisible. (Here we are using infinitely divisible as a synonym for associated.)

\def\noopsort#1{} \def\printfirst#1#2{#1}
\def\singleletter#1{#1}
            \def\switchargs#1#2{#2#1}
\def\bibsameauth{\leavevmode\vrule height .1ex
            depth 0pt width 2.3em\relax\,}
\makeatletter
\renewcommand{\@biblabel}[1]{\hfill#1.}\makeatother
\newcommand{\bysame}{\leavevmode\hbox to3em{\hrulefill}\,}

\def\wh{\widehat}
\def\ol{\overline}
\end{document}